\documentclass{article}

\usepackage{arxiv}

\usepackage[utf8]{inputenc} 
\usepackage[T1]{fontenc}    
\usepackage{hyperref}       
\usepackage{url}            
\usepackage{booktabs}       
\usepackage{amsfonts}       
\usepackage{nicefrac}       
\usepackage{microtype}      
\usepackage{lipsum}
\usepackage{graphicx}
\graphicspath{ {./images/} }

\title{On solvability of the non-local problem for the  fractional mixed-type equation with Bessel operator}

\author{Bakhodirjon Toshtemirov \\
  V.I.Romanovskiy Institute of Mathematics, Uzbekistan and \\
  Ghent University,  Belgium\\
  \texttt{bakhodirjon.toshtemirov@ugent.be}} \\

\begin{document}
\maketitle
\begin{abstract}
The non-local problem is considered for the partial differential equation of mixed-type with Bessel operator and fractional order. Explicit solution is represented by Fourier-Bessel series in the given domain. It is established the connection between the given data and the uniquely solvability of the problem. 
\end{abstract}

\keywords{Mixed-type equation \and the bi-ordinal Hilfer derivative \and hyper-Bessel fractional differential opertor}

\section{Introduction}
Fractional differential equations plays a significant role, because of its multiple applications in  engineering, chemistry, biology and other parts of science and modeling the real-life problems \cite{1}-\cite{3}. Studying boundary value problems for linear and non-linear fractional differential equations with Riemann-Liouville and Caputo fractional derivative have becoming interesting targets simultaneously \cite{4}-\cite{7}.

The quality and the types of articles have been changed when the generalized Riemann-Liouville differential operators (later called Hilfer derivative) used in the scientific field with its interesting application \cite{8}, \cite{9}. To tell the truth, the generalization of the Riemann-Lioville fractional differential operators was already announced by M. M. Dzherbashian and A. B.
Nersesian \cite{10} in 1968, but because of some reasons it was not familiar with many mathematicians around the globe till the translation of this work published in the journal of Fract. Calc. Appl.
An. \cite{11}.
 We also refer some papers \cite{12}-\cite{13} devoted studying some problems with the Dzherbashian-Nersesian differential operator which has the following form
\begin{equation}\label{1}
D_{0x}^{\sigma_n}g(x)=I_{0x}^{1-\gamma_n}D_{0x}^{\gamma_{n-1}}...D_{0x}^{\gamma_{1}}D_{0x}^{\gamma_{0}}, ~ ~ n\in\mathbb{N}, ~ x>0
\end{equation}
where $I_{0x}^{\alpha}$ and $D_{0x}^{\alpha}$ are the Riemann-Liouville fractional
integral and the Riemann-Liouville fractional derivative of order $\alpha$ respectively,  $\sigma_n\in(0, n]$ which is defined by
$$
\sigma_n=\sum\limits_{j=0}^{n}\gamma_j-1>0, ~ \gamma_j\in(0, 1].
$$

 V.I.Bulavatsky considered specific generalization of Hilfer fractional derivative \cite{14}, the particular case of the Dzherbashian-Nersesian differential operator which forms as follows \cite{15}:
\begin{equation}\label{2}
	D_{0\pm}^{(\alpha, \beta)\mu}g(t)=
I_{0\pm}^{\mu(n-\alpha)}\left(\pm\frac{d}{dt}\right)^nI_{0\pm}^{(1-\mu)(n-\beta)}g(t)
\end{equation}
This operator is similar to Hilfer derivative in terms of its interpolation concept between the Riemann-Liouville and Caputo fractional derivatives, specifically, i.e.
$$
D_{a+}^{(\alpha, \beta)\mu}g(t)=\left\{ \begin{gathered}
	D_{a+}^{\beta}g(t), \, \, \, \, \mu=0,\hfill \cr
	^CD_{a+}^{\alpha}g(t), \, \, \, \, \, \mu=1. \hfill \cr
\end{gathered}  \right.
$$

In \cite{15} the Cauchy problem is investigated for the ordinary differential equation involving the right-sided bi-ordinal Hilfer fractional derivative:
        \begin{equation}\label{3}
	\left\{\begin{array}{l} D_{0-}^{(\alpha, \beta)\mu}u(t)=\lambda u(t)+g(t),\\ \lim\limits_{t\to 0-}I_{0-}^{2-\gamma}u(t)=\xi_0,\\ \lim\limits_{t\to 0-}\frac{d}{dt}I_{0-}^{2-\gamma}u(t)=\xi_1,
	\end{array}\right.
\end{equation}
where $1<\alpha, \beta\leq2$, $\gamma=\beta+\mu(2-\beta)$, $\xi_0, \xi_1\in \mathbb{R},$ $g(t)$ is the given function.

\textbf{Lemma 1.}
Let $g(t)\in C[-T, 0], ~ g'(t)\in L_{1}(-T, 0)$. Then the solution of the problem (\ref{3}) as follows:
$$
u(t)=\xi_0(-t)^{\gamma-2}E_{\delta, \gamma-1}[\lambda(-t)^{\delta}]-\xi_1(-t)^{\gamma-1}E_{\delta, \gamma}[\lambda(-t)^{\delta}]+
$$
\begin{equation}\label{4}
	+\int\limits_{t}^{0}(z-t)^{\delta-1}E_{\delta, \delta}\left[\lambda(z-t)^{\delta}\right]g(z)dz,
\end{equation}
where $\delta=\beta+\mu(\alpha-\beta)$,  $\gamma=\beta+\mu(2-\beta)$.

We notice that while have been investigated the initial-boundary value or non-local problems involving popular class of differential operators like Riemann-Liouville, Caputo \cite{16}, \cite{17}, Hilfer fractional derivatives,  Hadamard \cite{18}, Hilfer-Hadamard, Prabhakar, Atangana-Baleanu \cite{19} the interest to another type of the differential operators is increased by many scientists, for instance, the hyper-Bessel differential operator \cite{20}, \cite{21} is becoming main target of research.
The importance of the hyper-Bessel differential operator is increasing since introduced by Dimovski \cite{22} because of its applications in science.  For example, in \cite{23} authors used hyper-Bessel differential operator to investigate  heat diffusion equation for discribing Brownian motion. In \cite{24} it is investigated fractional relaxation equation the regularized Caputo-like counterpart of the hyper-Bessel operator which has the following form
\begin{equation}\label{5}
^{C}\Big(t^{\theta}\frac{d}{dt}\Big)^\alpha{f}(t)=(1-\theta)^\alpha{t}^{-\alpha(1-\theta)}D_{1-\theta}^{-\alpha, \alpha}\left(f(t)-f(0)\right)
\end{equation}
where
 $$
D_{\beta}^{\gamma, \delta}f(t)\prod_{j=1}^{n}\left(\gamma+j+\frac{t}{\beta}\frac{d}{dt}\right)\big(I_{\beta}^{\gamma+\delta, n-\delta} f(t)\big),
$$
  is a  Erdelyi-Kober fractional derivative and $I^{\gamma, \delta}_{\beta}$ is a Erdelyi-Kober (E-K) fractional integral  defined as \cite{25}
  $$
I^{\gamma, \delta}_{\beta}f(t)=\frac{t^{-\beta(\gamma+\delta)}}{\Gamma(\delta)}\int_{0}^{t}(t^\beta-\tau^\beta)^{\delta-1}\tau^{\beta\gamma}f(\tau)d(\tau^\beta),
$$
which can be reduced up with weight to  $I^{q}_{0+}f(t)$ Riemann-Liouville fractional integral at $\gamma=0$ and $\beta=1$.
In \cite{26}  Fatma Al-Musalhi, et al. considered direct and inverse problems for the fractional diffusion equation with the regularized Caputo-like counterpart of the hyper-Bessel operator and proved the theorem of existence and uniqueness of the solution.

         Several local and nonlocal boundary value problems for mixed-type equations, ie, elliptic-hyperbolic and hyperbolic type equations were published \cite{27}-\cite{29}. The interesting point is that the conjugation conditions are taken according to the considered mixed-type equations and domains \cite{30}, \cite{31}.

\section{Formulation of the problem and its investigation}
\label{sec:headings}
In the present work, we investigate the following fractional differential equation of mixed type involving the regularized Caputo-like counterpart of the hyper-Bessel operator and the bi-ordinal Hilfer derivative in  $\Omega=\Omega_1\cup\Omega_2\cup{Q}$ domain:
  \begin{equation}\label{6}
 	f(x, t) = \left\{ \begin{gathered}
 		^C\left( {{t^\theta }\frac{\partial }{{\partial t}}} \right)^{\alpha_1} u\left( {x,t} \right) -\frac{1}{x}u_{x}(x, t)-
 		{u_{xx}}\left( {x,t} \right),\,\,\,\,(x,t)\in{\Omega_1},  \hfill \cr
 		{D^{(\alpha_2, \beta_2)\mu}_{0-}u}\left( {x,t} \right) -\frac{1}{x}u_{x}(x, t)- {u_{xx}}\left( {x,t} \right),\,\,\,\,(x,t)\in{\Omega_2},  \hfill \cr
 	\end{gathered}  \right.
 \end{equation}
where  $\Omega_1=\{(x,t): 0<x<1,\, \, 0<t<T\}$,  \, \,   $\Omega_2=\{(x,t): 0<x<1,\, -T<t<0\}$,  \, \,
 $Q=\{(x,t): 0<x<1, \,   t=0 \},$    \, \, $0<\alpha_1\leq{1}$, \, \, $\theta<1$, \,  \, $1<\alpha_2, \beta_2\leq2$, \, \,  $0\leq\mu\leq1$,\\
 $D^{(\alpha_2, \beta_2)\mu}_{0-}$ is the right-sided bi-ordinal Hilfer fractional derivative in the form (\ref{2}) and  \\
 $^{C}\Big(t^{\theta}\frac{d}{dt}\Big)^{\alpha_1}$ is the regularized Caputo-like counterpart of the hyper-Bessel fractional differential operator defined as (\ref{5}).

\textbf{Problem.}
Find a solution of eq.(\ref{6}) in $\Omega$, which satisfies the following regularity conditions
\begin{center}
$u(x,t)\in{C}(\Omega)$, \, \,  $u(\cdot, t)\in{C^1}\cap{C^2}(\Omega_1\cup\Omega_2)$, \\ $ ^C\Big(t^\theta\frac{\partial}{\partial{t}}\Big)^\alpha{u}(x,t)\in{C}(\Omega_1)$,
 \, \, \, $D_{0-}^{(\alpha_2, \beta_2)\mu}u(x,t)\in{C(\Omega_2)}$
\end{center}
along with the boundary conditions
\begin{equation}\label{7}
\lim_{t\to 0+}xu_{x}(x,t)=0, \, \, \, u(1,t)=0,  \, \,  \,   \, \,  \,   -T\leq{t}\leq{T},
\end{equation}
non-local condition
\begin{equation}\label{8}
\sum_{i=1}^{m}p_iI^{(1-\mu)(2-\beta_2)}_{0-}u(x,\xi_i)=u(x, T), \, \, \, \, 0\leq{x}\leq1,
\end{equation}
and the gluing conditions
\begin{equation}\label{9}
\mathop {\lim }\limits_{t \to  0-}{I^{(1-\mu)(2-\beta_2)}_{0-}}u(x,t)=\mathop {\lim }\limits_{t \to  0+}u(x,t), \, \, \, \, \,  0\leq{x}\leq{1},
\end{equation}

\begin{equation}\label{10}
\mathop {\lim }\limits_{t \to  0-}\frac{d}{dt}{I^{(1-\mu)(2-\beta_2)}_{0-}}u(x,t)=\mathop {\lim }\limits_{t \to  0+}t^{1-(1-\theta)\alpha_1}u_t(x,t), \, \, \, \, \,  0<x<1,
\end{equation}
where $-T\leq\xi_1<\xi_2<....<\xi_m\leq0, $ \, $f(x,t)$ is a given function.

 Using the separation variables method we obtain the following spectral problem
    \begin{equation}\label{11}
      X''(x)-\frac{1}{x}X'(x)+\lambda^2 X(x)=0,
    \end{equation}
        \begin{equation}\label{12}
          \lim\limits_{x\to 0+}xX(x)=0, X(1)=0,
        \end{equation}

  (9) is a Bessel equation of order zero; furthermore, (9), (10)  is a self-adjoint problem and its eigenfunctions are the Bessel functions given as follows
    \begin{equation}\label{13}
      X_k(x)=J_{0}(\lambda_k x), \, \, k=1, 2, ...,
    \end{equation}
    and the eigenvalues $\lambda_k$, are the positive zeros of $J_{0}(x)$,  i.e.
    $$
    \lambda_k=\pi k-\frac{\pi}{4}.
    $$
        The system of eigenfunctions $\{X_k\}$ forms an orthogonal basis in $L^2(0, 1)$, hence we can write sought function and given function in the form of series expansions as follows:
\begin{equation}\label{14}
  u(x, t)=\sum\limits_{k=1}^{\infty}u_k(t)J_{0}(x),
\end{equation}
  \begin{equation}\label{15}
    f(x, t)=\sum\limits_{k=1}^{\infty}f_k(t)J_{0}(x),
  \end{equation}
where $u_k(t)$ is not knownt yet, and $f_k(t)$ is the coefficient of Fourier-Bessel series, i.e.
$$
f_k(t)=\frac{2}{J_{1}^{2}(\lambda_k)}\int\limits_{0}^{1}x f(x, t)J_{0}(\lambda_k x)dx.
$$

Let us introduce new notations:
\begin{equation}\label{16}
\mathop {\lim }\limits_{t \to  0-}{I^{(1-\mu)(2-\beta_2)}_{0-}}u(x,t)=\varphi(x), \, \, \,  0\leq x\leq1,
\end{equation}
\begin{equation}\label{17}
\mathop {\lim }\limits_{t \to  0-}\frac{d}{dt}{I^{(1-\mu)(2-\beta_2)}_{0-}}u(x,t)=\psi(x), \, \, \,  0<x<1,
\end{equation}

\begin{equation}\label{18}
	\mathop {\lim }\limits_{t \to  0+}u(x,t)=\tau(x),  \, \, \, \,  0\leq{x}\leq1,
\end{equation}
here $\varphi(x), \, \, \, \tau(x)$ and $\psi(x)$ are unknown functions to be found later.

Further, after substituting  (\ref{14}) and (\ref{15}) into the eq.(\ref{6})  and initial conditions (\ref{16}), (\ref{17}), (\ref{18}), we obtain the following prdoblems
  \begin{equation}\label{19}
   \left\{\begin{gathered}
     ^C\Big(t^{\theta}\frac{d}{dt}\Big)^{\alpha_1}u_k(t)+\lambda_k^2 u_k(t)=f_k(t),
   \hfill \cr
u_k(0+)=\tau_k,
\end{gathered}\right.
\end{equation}
and
  \begin{equation}\label{20}
   \left\{\begin{gathered}
   D_{0-}^{(\alpha_2, \beta_2)\mu}u_k(t)+\lambda_k^2 u_k(t)=f_k(t),
   \hfill \cr
I_{0-}^{(1-\mu)(2-\beta_2)}u_k(0-)=\varphi_k,
\hfill \cr
\frac{d}{dt}I_{0-}^{(1-\mu)(2-\beta_2)}u_k(0-)=\psi_k,
\end{gathered}\right.
\end{equation}
in $\Omega_1$ and $\Omega_2$ respectively.

 The problem (\ref{19}) was studied in \cite{26} and we can write the solution of the considered problem in $\Omega_1$:
\begin{equation}\label{21}
u(x,t)=\sum_{k=1}^{\infty}\Big[\tau_kE_{\alpha_1, 1}\Big(-\frac{{\lambda_k}^2}{p^{\alpha_1}}t^{p\alpha_1}\Big)+G_k(t)\Big]J_{0}(\lambda_k x),
\end{equation}
here $p=1-\theta$ and
$$
G_k(t)=\frac{1}{p^{\alpha_1}\Gamma(\alpha_1)}\int_{0}^{t}\big(t^p-\tau^p\big)^{\alpha_1-1}f_k(\tau)d(\tau^p)-
$$
$$
-\frac{{\lambda_k}^2}{p^{2\alpha_1}}\int_{0}^{t}\big(t^p-\tau^p\big)^{2\alpha_1-1}E_{\alpha_1, 2\alpha_1}\Big[-\frac{{\lambda_k}^2}{p^{\alpha_1}}(t^p-\tau^p)^{\alpha_1}\Big]f_k(\tau)d(\tau^p),
$$
where $\tau_k$ is not known yet.

As we mentioned above the problem (\ref{20}) studied in \cite{15} and using the Lemma 1 we can write the solution of (\ref{20}) and then, considering (\ref{14}) the solution of considered problem in $\Omega_2$ domain  can be represented in the following form
\[
 u(x, t)=\sum\limits_{k=1}^{+\infty}\varphi_k(-t)^{\gamma_2-2}E_{\delta_2, \gamma_2-1}[-\lambda^2_k(-t)^{\delta_2}]J_{0}(\lambda_kx)-
 \]
 \begin{equation*}
 -\sum\limits_{k=1}^{+\infty}\psi_k(-t)^{\gamma_2-1}E_{\delta_2, \gamma_2}[-\lambda^2_k(-t)^{\delta_2}]J_{0}(\lambda_kx)+	
  \end{equation*}
 \begin{equation}\label{22}
 	+\sum\limits_{k=1}^{+\infty}\int\limits_{t}^{0}(z-t)^{\delta_2-1}E_{\delta_2, \delta_2}[-\lambda^2_k(z-t)^{\delta_2}]f_k(z)dzJ_{0}(\lambda_kx),
 \end{equation}
where $\gamma_2=\beta_2+\mu(2-\beta_2)$, $\delta_2=\beta_2+\mu(\alpha_2-\beta_2)$ and $\varphi_k$, $\psi_k$ are not known yet.

 After substituting (\ref{21}) and (\ref{22}) into gluing conditions with considering (\ref{14}), (\ref{15}) we obtain the following system of equations with respect to $\tau_k$, $\varphi_k$ and $\psi_k$:
\begin{equation}\label{23}
  \left\{\begin{array}{cc}
    \psi_k=-\frac{\lambda^2_k}{\Gamma(\alpha_1)}\tau_k, & \\ \\ \tau_k=\varphi_k
  \end{array}\right.
\end{equation}

Considering non-local condition (\ref{8}) and from (\ref{23}) we find unknowns as follows
\begin{equation}\label{24e}
  \tau_k=\varphi_k=\frac{F_k}{\Delta_k},
\end{equation}
\begin{equation}\label{25e}
  \psi_k=\frac{-\lambda_k^2}{p^{\alpha_1}\Gamma(\alpha_1)}\frac{F_k}{\Delta_k},
\end{equation}
 where
$$
\Delta_k=\sum_{i=1}^{m}p_i\Big[E_{\delta_2, 1}(-\lambda_k^2(-\xi_i))+\frac{\lambda_k^2\xi_i}{p^{\alpha_1}\Gamma(\alpha_1)}E_{\delta_2, 2}(-\lambda_k^2(-\xi_i))\Big]-E_{\alpha_1, 1}\Big(-\frac{\lambda_k^2}{p^{\alpha_1}}T^{\alpha_1p}\Big)
$$
$$
F_k=G_k(T)-\sum_{i=1}^{m}p_i\int\limits_{\xi_i}^{0}(s-\xi_i)^{\delta_2-\gamma_2+1}E_{\delta_2, \delta_2-\gamma_2+2}(-\lambda_k^2(s-\xi_i)^{\delta_2})f_k(s)ds,
$$
$$
G_k(T)=\frac{1}{p^{\alpha_1}\Gamma(\alpha_1)}\int_{0}^{T}\big(T^p-\tau^p\big)^{\alpha_1-1}f_k(\tau)d(\tau^p)-
$$
$$
-\frac{{\lambda_k}^2}{p^{2\alpha_1}}\int_{0}^{T}\big(T^p-\tau^p\big)^{2\alpha_1-1}E_{\alpha_1, 2\alpha_1}\Big[-\frac{{\lambda_k}^2}{p^{\alpha_1}}(T^p-\tau^p)^{\alpha_1}\Big]f_k(\tau)d(\tau^p).
$$

If $\Delta_k\neq0$, then we can find $\tau_k$, $\varphi_k$, $\psi_k$ unknowns uniquely.

First we show that $\Delta_k\neq0$ for sufficiently large $k$. For this intention we use the following lemma obtained from the properties of Wright-type function studied by A. Pskhu in \cite{32}.

\textbf{Lemma 2.} \cite{32} If $\pi\geq|argz|>\frac{\pi\alpha}{2}+\varepsilon, \, \, \, \varepsilon>0$,  then the following relations are valid for $z \to \infty$:
$$
\mathop {\lim }\limits_{\mid{z}\mid \to  \infty} E_{\alpha, \beta}(z)=0,
$$
$$
\mathop {\lim }\limits_{\mid{z}\mid \to  \infty} zE_{\alpha, \beta}(z)=-\frac{1}{\Gamma(\beta-\alpha)}.
$$

 By using the lemma 2,  we can calculate the behaviour of $\Delta_k$ at $k \to \infty$
$$
 \mathop {\lim }\limits_{k \to  +\infty}\Delta_k=\mathop {\lim }\limits_{\mid{z_1}\mid \to  +\infty}\sum\limits_{i=1}^{m}p_i\big[E_{\delta_2, 1}(z_1)+\frac{1}{\Gamma(\alpha_1)p^{\alpha_1}}z_1E_{\delta_2, 2}(z_1)\big]-
 $$
 $$
 -\mathop {\lim }\limits_{\mid{z_2}\mid \to  +\infty}E_{\alpha_1, 1}(z_2)=
 \sum\limits_{i=1}^{m}\frac{p_i}{\Gamma(\alpha_1)p^{\alpha_1}\Gamma(2-\delta_2)}
$$
where $z_1=-\lambda_k^2(-\xi_i)$, $z_2=-\frac{\lambda_k^2}{p^{\alpha_1}}T^{\alpha_1p}$, $\lambda_k=\pi k-\frac{\pi}{4}$.

From the last equality it is seen that $\Delta_k\neq{0}$ for sufficiently large $k$.

\subsection{Unequeness of the solution}
To show the uniqueness of the solution, it is enough to prove that homogeneous problem has a trivial solution.

Let us first consider the following integral
\begin{equation}\label{24}
u_k(t)=\frac{2}{J_{1}^{2}(\lambda_k)}\int\limits_0^1 xu(x, t)J_0(\lambda_k x)dx, \, k=1, 2, 3, ....,
\end{equation}

Then we introduce another function based on (\ref{24})
\begin{equation}\label{25}
v_{\varepsilon}(t)=\frac{2}{J_{1}^{2}(\lambda_k)}\int\limits_{\varepsilon}^{1-\varepsilon} xu(x, t)J_0(\lambda_k x)dx, \, k=1, 2, 3, ....,
\end{equation}

Applying $^C\Big(t^{\theta}\frac{\partial}{\partial t}\Big)^{\alpha_1}$ and $D_{0-}^{(\alpha_2, \beta_2)\mu}$ to (\ref{25}) and using the equation (\ref{6}) in homogeneous case with respect to $t$

$$
^C\Big(t^{\theta}\frac{\partial}{\partial t}\Big)^{\alpha_1}v_{\varepsilon}(t)=\frac{2}{J_{1}^{2}(\lambda_k)}\int\limits_{\varepsilon}^{1-\varepsilon} {^C\Big(t^{\theta}\frac{\partial}{\partial t}\Big)^{\alpha_1}u(x, t) x J_0(\lambda_k x)dx}=
$$
$$
=\frac{2}{J_{1}^{2}(\lambda_k)}\int\limits_{\varepsilon}^{1-\varepsilon}\Big[u_{xx}(x, t)+\frac{1}{x}u_x(x,t)\Big] x J_0(\lambda_k x)dx=\frac{-2\lambda_k^2}{J_{1}^{2}(\lambda_k)}\int\limits_{\epsilon}^{1-\epsilon}u(x, t) x J_{0}(\lambda_k x)dx
$$

$$
D_{0-}^{(\alpha_2, \beta_2)\mu}v_{\varepsilon}(t)=\frac{2}{J_{1}^{2}(\lambda_k)}\int\limits_{\varepsilon}^{1-\varepsilon}D_{0-}^{(\alpha_2, \beta_2)\mu}u(x, t) x J_0(\lambda_k x)dx=
$$
$$
=\frac{2}{J_{1}^{2}(\lambda_k)}\int\limits_{\varepsilon}^{1-\varepsilon}\Big[u_{xx}(x, t)+\frac{1}{x}u_x(x,t)\Big] x J_0(\lambda_k x)dx=\frac{-2\lambda_k^2}{J_{1}^{2}(\lambda_k)}\int\limits_{\epsilon}^{1-\epsilon}u(x, t) x J_{0}(\lambda_k x)dx
$$
and integrating by parts twice the right sides of the equalities on $t\in(0, T)$ and $t\in(-T, 0)$, respectively, and passing to the limit on $\varepsilon \to +0$ yield
\begin{equation}\label{26}
\left\{\begin{array}{l} ^C\Big(t^{\theta}\frac{d}{d t}\Big)^{\alpha_1}u_{k}(t)+\lambda^2u_k(t)=0,  \, \, \,   t>0, \\ D_{0-}^{(\alpha_2, \beta_2)\mu}u_{k}(t)+\lambda^2u_k(t)=0, \, \, \, t<0. \end{array}\right.
\end{equation}

Considering conditions in (\ref{19}), (\ref{20}) in homogeneous case, (\ref{26}) has a solution $u_k(t)=0$ if $\Delta_k\neq0$. Then from (\ref{24}) and the completeness of the system $X_k(x)$ in the space $L_2[0, l]$, $u(x, t)\equiv0$ in $\overline{\Omega}$. This completes the prove of uniqueness of the solution of the main problem.

\subsection{Existence of the solution}
In the below, we present the well-known lemma about the upper bound of the Mittag-Leffler function and the theorem related to Fourier-Bessel series for showing the existence of the solution.

\textbf{Lemma 3.} Let $\alpha<2, \, \beta\in{R}$ and $\frac{\pi\alpha}{2}<\mu<min\{\pi, \pi\alpha\}$, $M^*>0$. Then,the following estimate hold
$$
|E_{\alpha, \beta}(z)|\leq\frac{M^*}{1+|z|}, \, \, \, \mu\leq|argz|\leq\pi, \, \, |z|\geq0.
$$

\textbf{Theorem 1.}
    Let $f(x)$ be a function defined on the interval [0, 1] such that $f(x)$ is differentiable $2s$ times $(s\geq1)$ and
\begin{itemize}
  \item $f(0)=f'(0)=....=f^{(2s-1)}(0)=0$
  \item $f^{(2s)}(x)$ is bounded (this derivative may not exist at certain points)
  \item $f(1)=f'(1)=...=f^{(2s-2)}(1)=0$
\end{itemize}
then the following inequalities satisfied by the Fourier-Bessel coefficients of $f(x)$:
\begin{equation*}
  |f_k|\leq\frac{M}{\lambda_k^{2s-\frac{1}{2}}}.
\end{equation*}

Considering above lemma 3 and theorem 1, we show the upper bound of $G_k(t)$:
$$
\mid{G_k(t)}\mid\leq \int\limits_{0}^{t}|t^p-\tau^p|^{\alpha_1-1}|f_k(\tau)|d(\tau^p)+
$$
$$
+\frac{\lambda_k^2}{p^{2\alpha_1}}\int\limits_{0}^{t}|t^p-\tau^p|^{2\alpha_1-1}|E_{\alpha_1, 2\alpha_1}\big[-\frac{\lambda_k^2}{p^{\alpha_1}}(t^p-\tau^p)^{\alpha_1}\big]||f_k(\tau)|d(\tau^p)\leq
$$
$$
\leq\frac{1}{p^{\alpha_1}\Gamma(\alpha_1)}\int\limits_{0}^{t}|t^p-\tau^p|^{\alpha_1-1}\frac{M}{{\lambda_k}^{7/2}}d(\tau^p)
+\frac{\lambda_k^2}{p^{2\alpha_1}}\int\limits_{0}^{t}\frac{p^{\alpha_1}|t^{p}-\tau^{p}|^{2\alpha_1-1} M^*}{p^{\alpha_1}+\lambda_{k}^2|t^{p}-\tau^p|^{\alpha_1}}\frac{M}{{\lambda_k}^{7/2}}d(\tau^p)\leq
$$
$$
\leq\Big[\frac{M}{{\lambda_k}^{7/2} p^{\alpha_1}\Gamma(\alpha_1)}+\frac{M^* M}{p^{\alpha_1} {\lambda_k}^{7/2}}\Big]\int\limits_{0}^{t}|t^p-\tau^p|^{\alpha_1-1}d(\tau^p)\leq\frac{1}{{\lambda_k}^{7/2}}\Big[\frac{M}{ p^{\alpha_1}\Gamma(\alpha_1)}+\frac{M^* M}{p^{\alpha_1} }\Big]\frac{|t|^{\alpha_1p}}{\alpha_1}\leq
$$
$$
\leq\frac{M_1}{{\lambda_k}^{7/2}}|t|^{\alpha_1p}\leq\frac{M_1}{{\lambda_k}^{7/2}}|T|^{\alpha_1p}, \, \, ~~\, \, \,~ \, \,  \, M_1=\frac{M}{ p^{\alpha_1}\Gamma(\alpha_1+1)}+\frac{M^* M}{\alpha_1 p^{\alpha_1} },
$$

By using the last inequality, lemma 3 and theorem 1, we can write the upper bound of $F_k$:
$$
|F_k|\leq |G_k(T)|+\sum\limits_{i=1}^{m}p_i \int\limits_{\xi_i}^{0}|s-\xi_i|^{\delta_2-\gamma_2+1}|E_{\delta_2, \delta_2-\gamma_2+2}(-{\lambda_k}^2(s-\xi_i)^{\delta_{2}})||f_k(s)|ds\leq
$$
$$
\leq \frac{M_1 T^{\alpha_1 p}}{{\lambda_k}^{7/2}}+\sum\limits_{i=1}^{m}p_i\int\limits_{\xi_i}^{0}|s-\xi_i|^{\delta_2-\gamma_2+1}\frac{M^*}{1+\lambda_k^2|s-\xi_i|^{\delta_2}}\frac{M}{\lambda_k^{7/2}}ds\leq
$$
$$
\leq \frac{M_1T^{\alpha_1 p}}{{\lambda_k}^{7/2}}+\sum\limits_{i=1}^{m}p_i\int\limits_{\xi_i}^{0}|s-\xi_i|^{1-\gamma_2}\frac{MM^*}{\lambda_k^{11/2}}\leq\frac{M_1T^{\alpha_1 p}}{{\lambda_k}^{7/2}}+\sum\limits_{i=1}^{m}p_i\frac{(-\xi_i)^{2-\gamma_2}}{(2-\gamma_2)}\frac{MM^*}{\lambda_k^{11/2}}\leq,
$$
$$
\leq\frac{M_2}{\lambda_k^{7/2}}, \, \, \, \, M_2=M_1T^{\alpha_1 p}+\sum\limits_{i=1}^{m}\frac{p_i M M^*(-\xi_i)^{2-\gamma_2}}{(2-\gamma_2) {\lambda_k}^{2}},
$$
or
\begin{equation}\label{30e}
  |F_k|\leq\frac{M_2}{\lambda_k^{7/2}}.
\end{equation}

Now considering (\ref{30e}) and $\Delta_k\neq0$ then, we write upper bounds of $\tau_k$, $\varphi_k$, $\psi_k$ in (\ref{24e}) and (\ref{25e}).
\begin{equation}\label{30}
  |\tau_k|=|\varphi_k|\leq|\frac{1}{\Delta_k}||F_k|\leq \frac{M_2}{|\Delta_k|\lambda_k^{7/2}},
\end{equation}
\begin{equation}\label{31}
  |\psi_k|=\Big|\frac{-\lambda_k^2}{p^{\alpha_1}\Gamma(\alpha_1)}\Big|\Big|\frac{F_k}{\Delta_k}\Big|\leq\frac{M_2}{p^{\alpha_1}\Gamma(\alpha_1)|\Delta_k|\lambda_k^{3/2}}.
\end{equation}

For proving the existence of the solution, we need to show uniform convergence of series representations of  $u(x, t)$, $u_x(x, t)$, $u_{xx}(x,t)$,  $^C\left(t^\theta\frac{\partial}{\partial{t}}\right)^{\alpha}u(x, t)$ and $D_{0-}^{(\alpha_2, \beta_2)\mu}u(x, t)$ by using the solution (\ref{21}) and (\ref{22}) in $\Omega_1$ and $\Omega_2$ respectively.

    According to the last inequality and theorem 1, then we can present the existence of the solutions in both domains.

    $$
|u(x,t)|=\sum_{k=1}^{\infty}|u_k(t)||J_0(\lambda_kx)|\leq \sum_{k=1}^{\infty}|u_k(t)|\leq \sum_{k=1}^{\infty}\Big[|\tau_k||E_{\alpha_1, 1}(-\frac{\lambda_k^2}{p^{\alpha_1}}t^{\alpha_1 p})|+|G_k(t)||\Big]\leq
$$
$$
\leq\sum_{k=1}^{\infty}\Big(\frac{p^{\alpha_1}}{p^{\alpha_1}+{\lambda_k}^2|t^{p\alpha_1}|}\frac{M_2}{|\Delta_k|\lambda_k^{7/2}}+\frac{M_1 T^{\alpha_1 p}}{{\lambda_k}^{7/2}}\Big)
$$
One can shows that the series representation of $u(x, t)$ is bounded by convergent numerical series and by Weierstrass
M-test,  the series of $u(x, t)$ converges uniformly in $\Omega_1$.

Now we remind some properties of Bessel functions \cite{33}:
$J'_{0}(x)=-J_1(x)$;\\
$2J'_1(x)=J_{0}(x)-J_{2}(x)$ and asymptotic formula
$
|J_{\nu}(\lambda)kx)|\leq \frac{2A}{\sqrt{\lambda_kx})}, ~ ~ \, \, \nu>-1/2, \, \, ~ \\A=const.
$

It is not difficult to see that the series representation of $u_{xx}(x, t)$ is bigger than $u_{x}(x, t)$ hence it is enough to show the uniform convergency of $u_{xx}(x, t)$.
By using these properties we have
$$
|u_{xx}(x, t)|\leq\sum\limits_{k=1}^{\infty}|u_k(t)|\Big|\frac{d^2}{dx^2}J_{0}(\lambda_k x)\Big|=\sum\limits_{k=1}^{\infty}|u_k(t)|\frac{\lambda_k^2}{2}\Big|J_{2}(\lambda_k x)-J_0(\lambda_k x)\Big|
$$

And as a similar way of $u(x, t)$ we can show that
$$
|u_{xx}(x, t)|\leq\sum\limits_{k=1}^{\infty}\Big(\frac{p^{\alpha_1}}{p^{\alpha_1}+{\lambda_k}^2|t^{p\alpha_1}|}\frac{M_2}{|\Delta_k|\lambda_k^{7/2}}+\frac{M_1 T^{\alpha_1 p}}{{\lambda_k}^{7/2}}\Big)\frac{2A}{\sqrt{\lambda_kx})}.
$$

From the last inequality we can see that the series representation of $u_{xx}(x, t)$ is bounded by convergent series. According to Weierstrass
M-test,  the series of $u_{xx}(x, t)$ converges uniformly in $\Omega_1$.

The uniform convergency of $^C\left(t^\theta\frac{\partial}{\partial{t}}\right)^{\alpha}u(x, t)$ which is defined as
$$
^C\left(t^\theta\frac{\partial}{\partial{t}}\right)^{\alpha}u(x, t)=u_{xx}(x, t)-\frac{1}{x}u_x(x, t)+f(x, t)
$$
is similar to the way of showing convergency of the series representation of $u_{xx}(x, t)$.

In $\Omega_2$ domain it is enough to show the uniform convergency of $u_{xx}(x, t)$ which is bigger than other series. Hence the convergency of the series of $u(x, t)$, $u_x(x, t)$, $D_{0-}^{(\alpha_2, \beta_2)\mu}u(x, t)$ can be derived from the uniform convergency of $u_{xx}(x, t)$.

From theorem 1 and lemma 3, in $\Omega_2$ we can have
$$
|u_{xx}(x, t)|\leq\sum\limits_{k=1}^{\infty}|u_k(t)|\frac{\lambda_k^2}{2}\Big|J_{2}(\lambda_k x)-J_0(\lambda_k x)\Big|\leq\sum_{k=1}^{\infty}\lambda_k^2|\varphi_k||(-t)^{\gamma_2-2}||E_{\delta_2, \gamma_2-1}(-\lambda_k^2(-t)^{\delta_2})|+
$$
$$
+\sum_{k=1}^{\infty}\lambda_k^2|\psi_k||(-t)^{\gamma_2-1}||E_{\delta_2, \gamma_2}(-\lambda_k^2(-t)^{\delta_2})|+\sum\limits_{k=1}^{+\infty}\lambda_k^2\int\limits_{t}^{0}|z-t|^{\delta_2-1}|E_{\delta_2, \delta_2}[-\lambda^2_k(z-t)^{\delta_2}]||f_k(z)|dz\leq
$$
$$
\leq\sum_{k=1}^{\infty}\Big[\frac{\lambda_k^2 M_2}{|\Delta_k|\lambda_k^{7/2}}\frac{|(-t)^{\gamma_2-2}|M^*}{1+\lambda_k^2|(-t)^{\delta_2}|}+\frac{\lambda_k^2 M_2}{|\Delta_k|\lambda_k^{3/2}p^{\alpha_1}\Gamma(\alpha_1)}\frac{|(-t)^{\gamma_2-1}|M^*}{1+\lambda_k^2|(-t)^{\delta_2}|}\Big]+
$$
$$
+\sum\limits_{k=1}^{+\infty}\lambda_k^2\int\limits_{t}^{0}|z-t|^{\delta_2-1}\frac{M^*}{1+\lambda_k^2|(z-t)^{\delta_2}|}\frac{M}{\lambda_k^{7/2}}dz\leq
$$
$$
\leq\sum_{k=1}^{\infty}\frac{1}{\lambda_k^{3/2}}\Big[\frac{M_2 M^* T^{\gamma_2-2}}{|\Delta_k|(1+{\lambda_k}^{2}T^{\delta_2})}+\frac{M_2 M^* T^{\gamma_2-\delta_2-1}}{|\Delta_k|}+\frac{M M^* \ln(1+\lambda_k^2 T^{\delta_2})}{\delta_2 \lambda_k^{2}}\Big]\leq\sum_{k=1}^{\infty}\frac{M_3}{\lambda_k^{3/2}}
$$
where $\lim\limits_{\lambda_k \to \infty} \frac{\ln(1+\lambda_k^2 T^{\delta_2})}{\delta_2 \lambda_k^{2}}<\infty$ according to l'Hopital's rule. It can be seen that the series representation of $u_{xx}(x, t)$ is bounded by convergent numerical series and due to Weierstrass
M-test,  the series of $u_{xx}(x, t)$ converges uniformly in $\Omega_2$

Using (\ref{30}), (\ref{31}) and theorem 1, lemma 3, we can show the uniform convergence of $u(x, t)$, $u_{x}(x,t)$,  and $D_{0-}^{(\alpha_2, \beta_2)\mu}u(x, t)$ as a similar method used for $u_{xx}(x, t)$.
    \smallskip

    Finally we have proved the uniqueness and existence of the solution to the considered problem as stated in the following theorem.
   
   \textbf{Theorem 2.}  Let $\Delta_k\neq0$ and $f(x,t)$ satisfies the following conditions
\begin{itemize}
  \item $f(0, t)=f'_x(0, t)=....=f'''_x(0, t)=0$;
  \item $f(1, t)=f'_x(1, t)=f''_x(1, t)=0$;
  \item $\frac{\partial^4}{\partial x^4}f(x, t)$ is bounded,
\end{itemize}
then the considered problem's unique solution exists.
    
\section{Acknowledgments}
The author is very grateful to his supervisor prof. Erkinjon Karimov for his valuable suggestions and to the reviewers in terms of improving the quality of this research.
\bibliographystyle{unsrt}  


\end{document}